\documentclass[11pt]{article} 
\usepackage{amsmath,amssymb, latexsym, amsthm, amscd}
\def\A{{\mathbb{A}}}
\def\C{{\mathbb{C}}}
\def\P{{\mathbb{P}}}
\def\Z{{\mathbb{Z}}}
\def\R{{\mathbb{R}}}
\def\H{{\mathbb{H}}}
\def\Hbar{{\overline{\mathbb{H}}}}
\def\Q{{\mathbb{Q}}}
\def\CC{{\mathcal{C}}}
\def\M{{\mathcal{M}}}

\def\O{{\mathcal{O}}}
\def\S{{\mathcal{S}}}
\def\T{{\mathcal{T}}}
\def\D{{\Delta}}
\def\a{{\alpha}}

\def\d{{\delta}}
\def\G{{\Gamma}}
\def\L{{\Lambda}}
\def\Ldual{{\Lambda^*}}
\def\l{{\lambda}}
\def\g{{\gamma}}
\def\r{{\rho}}
\def\s{{\sigma}}
\def\z{{\zeta}}
\def\w{{\omega}}
\def\b{{\backslash}}
\def\til{\widetilde}
\def\<{{\langle}}
\def\>{{\rangle}}
\def\isom{{\;\cong\;}}

\DeclareMathOperator{\tr}{tr}
\def\SL{{\rm SL}}
\def\slz{{{\rm SL}_2(\Z)}}
\def\slzn{{{\rm SL}_2(\Z/N)}}
\newcommand{\pslz}{\mathrm{PSL}_2(\Z)}
\def\abcd{{a\;\;b \choose c\;\;d}}
\newcommand{\twobytwo}[4]
{\left[\begin{array}{rr} #1 & #2 \\ #3 & #4 \end{array}\right]}

\newtheorem{theorem}{Theorem}

\newtheorem{prop}[theorem]{Proposition}

\begin{document}

\title{Level structures on the Weierstrass family of cubics}
\author{Mira Bernstein\\Department of Mathematics\\Wellesley College,\\Wellesley, Massachusetts, USA \and 
Christopher Tuffley\\Department of Mathematics\\University of California at Davis\\ Davis, California, USA}
\maketitle

\renewcommand{\thefootnote}{}

\begin{abstract}
Let $W\to\A^2$ be the universal Weierstrass family of cubic curves
over $\C$. For each $N\geq 2$, we construct surfaces parameterizing
the three standard kinds of level $N$ structures on the smooth fibers
of $W$. We then complete these surfaces to finite covers of
$\A^2$. Since $W\to \A^2$ is the versal deformation space of a cusp
singularity, these surfaces convey information about the level
structure on any family of curves of genus $g$ degenerating to a
cuspidal curve. Our goal in this note is to determine for which values
of $N$ these surfaces are smooth over $(0,0)$. From a topological
perspective, the results determine the homeomorphism type of 
certain branched covers of $S^3$ with monodromy in
$\slzn$.
\end{abstract}

\section{Introduction} \label{intro}
Let $W$ be the Weierstrass family of  plane cubic curves
$y^2=x^3+ax+b$, with base $\A^2_{a,b}$. Let $\D$ be the discriminant
locus in $\A^2$, defined by $4a^3+27b^2=0$. Over $\A^2\b\D$, $W$ has smooth
fibers and thus restricts to a family $\til{W}$ of elliptic curves. Working over $\C$, for each positive integer $N$ we construct the following finite \'{e}tale covers of $\A^2\b\D$: 
\begin{eqnarray}
\til{S}(N)&=&
 \{(p,q)\in \til{W}[N]\times_{\A^2\b\D} \til{W}[N]: \<p,q\>=e^{2\pi i/N}\}; 
   \label{SofN.eq} \\
\til{S}_1(N)&=& \{p\in\til{W}[N]: p\text{ has exact order }N\}; 
     \label{S1ofN.eq}\\
\til{S}_0(N)&=& 
  \{\gamma\subset\til{W}[N]:\gamma\text{ a cyclic subgroup of order }N\}.
  \label{S0ofN.eq}
\end{eqnarray}
These surfaces parametrize the three standard kinds of level $N$ structures on the smooth fibers of $W$.

What happens to the level structure as the smooth fibers degenerate to
singular ones -- in particular, to the cuspidal curve $y^2=x^3$? One
natural way to interpret this question is to take the normalization of
$\A^2$ in the function fields of
$\til{S}(N)$, $\til{S}_1(N)$, and $\til{S}_0(N)$ and to study the
resulting surfaces, which we denote by $S(N)$, $S_1(N)$ and
$S_0(N)$. We will use the shorthand $S_i(N)$ when making an assertion
that is valid for any one of the three spaces.

$S_i(N)$ is a normal finite cover of $\A^2$, ramified only over $\D$.  Let $U$ be a neighborhood of the origin in $\A^2$, and let $p\in U\b\D$. $U\b\D$ is homotopic to the complement of a trefoil knot in $S^3$, and it is not hard to see that the monodromy acts transitively on the points lying over $p$ in $S_i(N)$. Thus the fiber over the origin in each of our surfaces consists of a single point $Q$.  Our goal in this note is to prove the following

{\bf Theorem.} {\em For all $N$, $S_i(N)$ is smooth away from $Q$. At $Q$,
\begin{itemize}
\item[] $S(N)$ is smooth iff $N=2$ or $3$; 
\item[]$S_1(N)$ is smooth iff $\;2\leq N \leq 6$; 
\item[]$S_0(N)$ is smooth iff $N=2$ or $4$.
\end{itemize}}

The method of proof is as follows: for each $S_i(N)$, we first
construct a surface $T_i(N)$ with known singularities which has a
smooth open subset $U\isom S_i(N)\b Q$. We then determine whether the
divisor $T_i(N)\b U$ can be blown down to obtain a smooth surface. The
construction of $T_i(N)$ is described in Section 2, and the proof of
the Theorem is given in Sections 3 and 4.

There are several reasons for studying the surfaces $S_i(N)$, and
specifically their geometry at $Q$. One motivation is that $W\to \A^2$
is the versal deformation space of a cusp singularity (\cite{H-M},
p. 98). In other words, if $X\to U$ is the versal deformation space of
any cuspidal curve $C$ of genus $g$, then $U$ fibers over $\A^2$ and,
in the neighborhood of the cusp, $X$ fibers over $W$. For a given type
of level $N$ structure (which, for a curve of genus $g$, refers to
$N$-torsion data on the Jacobian), we can construct a finite cover
$\til{\S}(N)$ of $U\b\D$ that parametrizes structures of this type on
the smooth fibers of $X$, and then complete $\til{\S}(N)$ to a cover
$\S(N)$ of $U$. In general, $\S(N)$ will be reducible, and each of its
irreducible components will naturally fiber either over $\A^2$ or over
one of our surfaces $S_i(N)$ (or over some analogue of $S_i(N)$ for a
more unusual level structure).  Thus if we wish to understand {\em in
general} how level $N$ structures on smooth curves degenerate to
cuspidal curves, the family $W\to \A^2$ is the place to start.

Another motivation comes from low-dimensional topology. Let $D$ be the
boundary of a ball around the origin in $\C^2$. As a real manifold,
$D\isom S^3$. The preimage of $D$ in $S_1(N)$ is a finite cover of $D$
branched over $D\cap\D$, i.e. a finite cover of $S^3$ branched over a
trefoil knot $K$. Denote this 3-manifold by $M_1(N)$. There is an
extensive literature about such covers when the monodromy around $K$
is $\Z/N$ -- in other words, when we are given a one-dimensional
representation of $\pi_1(S^3\b K)$. Our case is, in some sense, the
next simplest: the monodromy is $\SL_2(\Z/N)$, based on the most
natural two-dimensional representation of the trefoil knot group. Our
results show that $M_1(N)$ is the three-sphere for $2\leq N\leq 6$,
and a circle bundle over a Riemann surface for $N\geq 7$. There are
known formulas for both the genus of the surface and the Chern class
of the bundle, so we have completely determined the homeomorphism type
of $M_1(N)$. Section \ref{topological.sec} gives an interpretation of
our results from this topological perspective.

{\bf Acknowledgements:} We are grateful to Barry Mazur for helping
formulate the problem and to Jordan Ellenberg, Joe Harris, and Mike
Roth for many helpful conversations.

\section{Constructing $T_i(N)$}

We begin by constructing the auxiliary surface $T(N)$ for
$N>2$. Recall that in this case, the modular curve $X(N)$
parameterizing elliptic curves with full level $N$ structure is a fine
moduli space.

\begin{prop}
For $N>2$, let $f:\CC \to X(N)$ be the universal family over $X(N)$,
with relative dualizing sheaf $\w_f$. Let $\L = f_*(\omega_f)$ be the
Hodge bundle on $X(N)$, and let $T(N)$ be the total space of the dual
bundle $\Ldual$, with zero section $Z$. Then $T(N)$ is the blow-up of
$S(N)$ at $Q$, with exceptional divisor $Z$.
\end{prop}

\begin{proof}
A Weierstrass equation $y^2 = x^3+ax+b$ for an elliptic curve $E$ corresponds to the choice of a non-vanishing holomorphic differential $\omega = dx/y$ on $E$. Any other Weierstrass equation for $E$ is of the form
$$
y^2 = x^3 + a_u x + b_u,
$$ 
where $a_u = u^4a$, $b_u = u^6b$ for some $u\neq 0$. The corresponding differential is
$\w_u = u^{-1}\w$. Let $\eta_u\in H^0(E,K_E)^*$ be the functional that sends  $\w_u$ to 1. Note that 
$$
\lim_{u\to 0} (a_u,b_u) = (0,0), \quad\lim_{u\to 0} \w_u = \infty, \quad\text{hence }
\lim_{u\to 0}\eta_u = 0.
$$

More generally, for an irreducible curve $C$ of arithmetic genus 1 (smooth or singular), the choice of a Weierstrass equation corresponds to a choice of a marked point $0\in C$ and of a non-zero linear functional $\eta$ on sections of the dualizing sheaf $\w_C$. The base $\A^2_{a,b} \b 0$ of the Weierstrass family (with the cuspidal curve removed) is thus a moduli space for such triples $(C,0,\eta)$. 
Similarly, $\til{S}(N)$ is a moduli space for the data
$$
[E,(p,q),\eta],
$$
where $E$ is an elliptic curve, $(p,q)$ is a full level $N$ structure on $E$, and $\eta\in H^0(E,\w_E)^*, \eta\neq 0$. 

Note that the fiber of $\Ldual$ over a point $[E, (p,q)]$ of $X(N)$ is precisely $H^0(E,\w_E)^*$.
Thus we have a natural inclusion $\iota: \til{S}(N)\to T(N)$, which
fits into the commutative diagram
$$
\begin{CD}
\til{S}(N) 	@>\iota >> 		T(N)\b Z\\
@VVV				@VVV\\
\A^2\b\D 	@>>> 		\A^2\b 0.	
\end{CD}
$$
Since $T(N)$ is smooth (hence normal), we conclude that $S(N)\b Q \;\isom \; T(N)\b Z$, by definition of $S(N)$. A neighborhood of $Q$ in $S(N)$ maps to a neighborhood of $Z$ in $T(N)$, so $T(N)$ is the blow-up of $S(N)$ at $Q$, as required.
\end{proof}

For $N>4$, we can define $T_1(N)$ in a similar way. However, this construction does not work for those $S_i(N)$ for which the corresponding modular curves are not fine moduli spaces. In that case, there is no universal family and hence no $\L$. Instead, we note that the natural $\slz$ action on $\til{S}(N)$ as a cover of $\A^2\b \D$ extends to $S(N)$, and hence to $T(N)$. We can then define $T_i(N)$ as a scheme-theoretic quotient: 
\begin{eqnarray*}
T(2) &=& T(4)/\G(2),\\
T_1(N)&=& T(N)/\G_1(N),\\
T_0(N)&=& T(N)/\G_0(N),
\end{eqnarray*}
where $\G(N)$, $\G_1(N)$, and $\G_0(N)$ are the usual congruence subgroups of level $N$.
In each case, let $Z'$ be the (reduced) image of the zero section from $T(N)$. Then one easily checks that
$$
T_i(N)\b Z' 
\isom (S(N)\b Q)/\G_i(N) \isom S_i(N)\b Q,
$$
and similarly for $T(2)$. Note that these quotient spaces may be singular due to isolated fixed points of the group action. To analyze the singularities, we will need to describe the $\slz$ action on $T(N)$ in terms of suitably chosen local coordinates.

Recall that $X(N) \isom \overline{\H}/\G(N)$, where 
$\overline{\H} = \H \cup \P^1(\Q)$ is the extended upper half plane. 
We would like to establish a similar model for $T(N)$.
Define $\psi: \H\times\C \to T(N)$ by
$$
\psi(\tau, t) = \left[E_\tau, \left(\frac{\tau}N, \frac1N\right), \;t\; dz^*\right],
$$
where $E_\tau=\C/\<1,\tau\>$ and $z$ is the coordinate on $E_\tau$ inherited from $\C$. 
Let $A=\abcd\in\slz$ act on $\H\times\C$ by
$$
A(\tau,t) = \left(\frac{a\tau+b}{c\tau+d},\frac{t}{c\tau+d}\right).
$$
It is not hard to see that $\psi(\tau,t)=\psi(\tau',t')$ if and only if $(\tau',t')$ and $(\tau,t)$ are in the same $\G(N)$-orbit under this action. Thus $\psi$ induces an inclusion of $(\H\times\C)/\G(N)$ as an open subset $\til{T}(N)$ of $T(N)$. 
The complement of $\til{T}(N)$ consists of fibers over the cusps of $X(N)$, which we will refer to as the {\em cusp fibers}. 

$\til{T}(N)$ has an $\slz$ action inherited from $\H\times\C$. If $A=\abcd\in\slz$ then
$A$ acts by
\begin{eqnarray*}
A\left[E_\tau,\;\left(\frac{\tau}{N},\frac1N\right),\;t\;dz^*\right] &=& 
\left[E_{A\tau},\;\left(\frac{A\tau}N,\frac1{N}\right),\;t\frac{dz^*}{c\tau+d}\right]\\
&=& \left[E_{\tau},\;\left(\frac{a\tau+b}N,\frac{c\tau+d}{N}\right),\;t\;dz^*\right].
\end{eqnarray*}
We will often write $\left(\frac{a\tau+b}N,\frac{c\tau+d}{N}\right)$ as $A{\tau/N \choose 1/N}$. Thus $A$ permutes the level structure without changing the differential. This is exactly the action of $A$ on $\til{S}(N)$ as a cover of $\A^2\b \D$. 

To see how this action extends to the cusp fibers of $T(N)$, we need to define local coordinates in the neighborhood of each cusp. Let $U$ be a disk of radius $<1$ and let $\z = e^{2\pi i/N}$. For each $M={k\;\;l \choose m\;\;n}\in\slz$, define a map $\psi_M: U\times\C \to T(N)$ by
$$
\psi_M(q,s) = [\C^*/\<q^N\>, \;(q^k\z^l, q^m\z^n), \; s\cdot 2\pi i u\;du^*],
$$
where $u$ is the coordinate on $\C^*/\<q^N\>$ inherited from $\C^*$. Note that when $q=0$, the differential $\frac{du}{2\pi i u}$ on $\C^*$ determines a non-vanishing section of the dualizing sheaf of $\P^1(\C)/_{0\sim\infty}$ and $2\pi u\;du^*$ is its dual, as required.  

If $q = e^{2\pi i \tau/N}$ then $u=e^{2\pi i z}$ gives an isomorphism
$E_\tau \to \C^*/\<q^N\>$. Thus
\begin{eqnarray*}
\psi_M(q,s) &=& \left[E_\tau, \left(\frac{k\tau+l}{N}, \frac{m\tau+n}{N}\right), \;s \;dz^*\right]\\
	&=& \left[E_{M\tau}, \left(\frac{M\tau}N, \frac1N\right), \;s\frac{dz^*}{m\tau+n}\right]\\
	&=& \psi\left(M\tau, \frac{s}{m\tau+n}\right).
\end{eqnarray*}
For a given value of $q\neq 0$, $\tau$ is only determined up to translation by $N$, but 
this does not affect the resulting point on $T(N)$. As $q$ goes to $0$, $\tau\to\infty$ and hence $M\tau\to\frac{k}{m}$. Thus $\psi_M$ gives a system of local coordinates on $T(N)$ near the cusp  $\frac{k}{m}$ on $Z\isom X(N)$. We denote these coordinates by $(q,s)_M$. We leave it to the reader to show that $A$ takes the point $(q,s)_M$ to the point $(q,s)_{AM}$; this is true for $q\neq 0$ and hence  for $q=0$ as well.

We are particularly interested in the fixed points of $T(N)$ under this action of $\slz$. The following two propositions deal with fixed points in $\til{T}(N)$ and on cusp fibers respectively:

\begin{prop} \label{fp1}
Suppose $\a = [E_\sigma, \left(\frac{\sigma}N,\frac1N\right), t\;dz^*]\in \til{T}(N)$ is fixed by $A\in\slz$. Then: 
\begin{enumerate} 
\item[(i)] $t=0$ (i.e. $\a\in Z$).
\item[(ii)] Let $\mathcal{T}_\a$ be the tangent space to $T(N)$ at $\a$. Let $x$ be the coordinate on $\T_\a$ corresponding to motion along $Z$ in $T(N)$, and let $t$ be the coordinate corresponding to motion along the fiber. Then there exists $\l\in\C$ such that multiplication by $\l$ induces an automorphism of $E_\sigma$ and the action of $A$ on $\T_\a$ is given by
$$
A(x,t) = \left(\frac{x}{\l^2}, \frac{t}{\l}\right).
$$
\end{enumerate}
\end{prop}

\begin{proof}
First, we may assume that $A\sigma = \s$ in $\H$. If not, then $A\s = \gamma\s$ for some $\g\in\G(N)$ and we can replace $A$ by $\g^{-1}A$ in the statement of the proposition. The action on $T(N)$ will be exactly the same.

If $A=\abcd$, let $\l=c\s+d$. Clearly multiplication by $\l$ gives an automorphism of $E_\s$, and
\begin{eqnarray*}
A\a&=&\left[E_{A\sigma}, \left(\frac{A\s}N, \frac{1}{N}\right), 
\frac{t}{\l}dz^*\right] \\
&=&\left[E_\sigma, \left(\frac{\s}N, \frac{1}{N}\right), \frac{t}{\l}\;dz^*\right]. 
\end{eqnarray*}
Comparing this with $\a$, we see that $\frac{t}{\l} = t$, so $t=0$, proving $(i)$. 

Let $U\times \C$ be a neighborhood of $(\s,0)$ in $\H\times\C$, small enough to be isomorphic onto its image in $T(N)$. The coordinates $(\tau, t)$ on $U\times \C$
define local coordinates $x=\tau-\s$ and $t$ on $T(N)$, centered at $\a$. Simple algebra shows that
$$
A(\s+x,t) = \left(\s+\frac x{\l^2(1+\frac c\l x)}, \;\frac t{\l(1+\frac{c}{\l}x)}\right).
$$

Quotienting by $\<x^2\>$ to pass to the tangent space, we obtain the stated action on $\T_\a$.
\end{proof}

\begin{prop} \label{fp2}
Let $\a\in Z$ be a cusp corresponding to the point $r\in\P^1(\Q)$. Choose any $M\in\slz$ such that $M\infty = r$ and define local coordinates $(q,s)_M$ at $\a$. If $A\in\slz$ fixes any point on the fiber over $\a$, then locally
$$
A(q,s) = (\z^k q, s) \quad\text{or}\quad A(q,s) = (\z^k q, -s),
$$
where $k\in \Z$ and $\z = e^{2\pi i/N}$. Thus $A$ fixes either the whole fiber or only $\a$ itself.
\end{prop}

\begin{proof}
As before, we can assume that $Ar=r$; otherwise $Ar=\g r$ for some $\g\in\G(N)$, and we can replace $A$ by $\g^{-1}A$. 

Let $B=M^{-1}AM$. Then $B\infty = \infty$, so 
$B=\pm{1\;\;k \choose 0\;\;1}$ for some $k\in Z$. If $q=e^{2\pi i \tau/N}$, then
\begin{eqnarray*}
A(q,s) 
	&=& \left[E_\tau, AM{\tau/N \choose 1/N}, s \;dz^*\right]\\
	&=& \left[E_\tau, MB{\tau/N \choose 1/N}, s\;dz^*\right]\\
	&=& \left[E_{\tau+k}, M{(\tau+k)/N \choose 1/N}, \pm s \; dz^*\right]\\
	&=&(\z^k q, \pm s).
\end{eqnarray*}

\end{proof}

Propositions \ref{fp1} and \ref{fp2} will be our main tools in analyzing the singularities 
of $T_i(N)$ in the next two sections.

\section{Proof of the Main Theorem for $S(N)$ and $S_1(N)$} \label{SN}

We first show that $S_i(N)$ is always smooth away from $Q$; this part of the proof is valid for $S_0(N)$ as well. For $N>2$, $T(N)$ is smooth, since it is the total space of a line bundle; hence $S(N)\b Q \isom T(N)\b Z$ is smooth as well. 
$T_1(N)$ and $T_0(N)$ may have quotient singularities due to isolated fixed points of the action of $\G_i(N)$. However, Propositions \ref{fp1} and \ref{fp2} show that isolated fixed points in $T(N)$ can occur only on $Z$. The corresponding singularities will be on the image $Z'$ of $Z$ in $T_i(N)$. Thus $S_i(N)\b Q \isom T_i(N)\b Z'$ is again smooth. The same argument works for $T(2)$ and $T_1(2)$ as quotients of $T(4)$.

Now consider smoothness at $Q$. For $N>2$, $S(N)$ is smooth at $Q$ if the zero section $Z$ can be blown down to obtain a smooth surface. For this to happen, $Z$ must be a rational curve with $Z^2 = -1$. In other words, $X(N)$ must have genus $0$ (so $N \leq 5$; see \cite[ch. 4] {Sch} for genus formulas for modular curves), and we must have $\L = \O(1)$. To compute the degree of $\L$, recall that on any one-parameter family of stable curves $f: Y\to X$, with relative dualizing sheaf $\w_f$, we have the identity
$$
12\l = \kappa+\d
$$
(\cite{H-M}, p. 158). Here $\l = c_1(\L)$, $\kappa = f_*(c_1(\w_f)^2)$ , and $\delta$ is the pullback to $X$ of the boundary divisor class $\d$ on the moduli space of stable curves. In the case of the universal family $\CC\to X(N)$, we have $\kappa = 0$ (since the restriction of $\w$ to every fiber is trivial) and $\d$ is the pullback of the point $j=\infty$ under the natural map $X(N)\to\overline{\M}_{1,1}$. The degree of $\d$ on $X(N)$ is the degree of $X(N)$ as a cover of the $j$-line, which is 
\begin{equation}
[\pslz:\Gamma(N)] = \frac12 {N^3}\prod_{p|N}\left(1-\frac1{p^2}\right).
\label{c1MN.eq}
\end{equation}
Thus the degree of $\l$ is $\frac{1}{12}[\pslz:\Gamma(N)]$ and $\L = \O(1)$ if and only if $N=3$. We have shown that $S(N)$ is smooth for $N=3$ and singular for $N>3$. 

The argument for $S_1(N)$, $N\geq 5$, is completely analogous. $X_1(N)$ is a fine moduli space, so we let $\L$ be the Hodge bundle for the universal family over $X_1(N)$, and let $T_1(N)$ be the total space of $\Ldual$. Then $T_1(N)$ is the blow-up of $S_1(N)$ at $Q$, with exceptional divisor $Z$, the zero section. For $S_1(N)$ to be smooth, $X_1(N)$ must have genus $0$, and we must have $\L = \O(1)$. The degree of $\l$ is now $\frac{1}{12}$ of the degree of $X_1(N)$ as a cover of the $j$-line, i.e.
\begin{equation}
\deg(\l) = \frac1{12}[\pslz: \Gamma_1(N)] = \frac1{24} {N^2}\prod_{p|N}\left(1-\frac1{p^2}\right).
\label{c1M1N.eq1}
\end{equation}
Thus $S_1(N)$ is smooth for $N=5$ and $6$ and singular for $N>6$. Note that the genus of $X_1(N)$ is not an obstruction to smoothness here or below, since $X_1(N)$ is rational for $N\leq 10$.

For $S_1(3)$ and $S_1(4)$, we must first determine what, if any, quotient singularities arise in $T_1(3)$ and $T_1(4)$. These singularities come from isolated fixed points of the $\G_1(N)$ action on $T(N)$, which all lie on $Z \isom X(N)$. Thus they correspond to isolated fixed points of the $\G_1(N)$ action on $\Hbar$, up to $\G_1(N)$ equivalence. 

$\G_1(3)$ has two such fixed points: $\infty$ and $\frac{1+e^{i\pi/3}}3$. Let $A={1\;\;1\choose 0\;\;1}$ be a generator for $\G_1(3)/\G(3)$. Setting $M=I$ in Proposition \ref{fp2}, we see that $A$ must fix the whole fiber over the cusp at $\infty$. Thus no singularities arise there. Now consider the action of $A$ near the point $\a = \psi(\frac{1+e^{i\pi/3}}3,0)\in Z$. Since $A$ has order 3, in Proposition \ref{fp1} we must have $\l^3=1$. The invariants of the action on the tangent space $\T_\a$ are $\{x^3, xy, y^3\}$. Thus the quotient $T_1(3)$ has a single $A_2$ singularity.

Now let $\pi: T(3)\to T_1(3)$ be the quotient map, and let $Z' = \pi(Z)$. Recall that $Z^2 = -1$ in $T(3)$. Then we have
\begin{eqnarray*}
&&\pi_*Z = 3Z',\; \pi^*Z' = Z, \\
&&\pi^*(Z'^2) = (\pi^*Z')^2 = Z^2 = -1,\\
&&Z'^2=-\frac13.
\end{eqnarray*}
Blowing up the $A_2$ singularity on $Z'$, we get exceptional divisors $E_1$ and $E_2$, each with self-intersection $-2$. The proper transform $\til{Z'}$ of $Z'$ has self-intersection $-1$. Thus we can successively blow down $\til{Z'}$, $E_1$ and $E_2$, obtaining a smooth space at each stage. We conclude that  $S_1(3)$ is  smooth.

We now repeat the same calculation for $T_1(4)$. $\G_1(4)$ has no fixed points in $\H$ and two (equivalence classes of) fixed points in $\P^1(\Q)$: $\frac 12$ and $\infty$. Once again, the generator $A={1\;\;1\choose 0\;\;1}$ of $\G_1(4)/\G(4)$ fixes the entire fiber over $\infty$. For the cusp at $\frac12$, set $M={1\;\;0\choose 2\;\;1}$ in Proposition \ref{fp2}. Then $MAM^{-1} = -{1 \;\;-1 \choose 0\;\;\;\;1}$ mod 4, so the action of $A$ is locally given by $(q,s)\to(-iq, -s)$. The invariants are $\{q^4, q^2s, s^2\}$, so the quotient $T_1(4)$ has a single $A_1$ singularity.

Once again, let $\pi: T(4)\to T_1(4)$ be the quotient map, and let $Z' = \pi(Z)$. Recall that $Z^2 = -2$ in $T(4)$. Then we have:
\begin{eqnarray*}
&&\pi_*Z = 4Z',\; \pi^*Z' = Z, \\
&&\pi^*(Z'^2) = (\pi^*Z)^2 = Z^2 = -2,\\
&&Z'^2=-\frac12.
\end{eqnarray*}
Blowing up the $A_1$ singularity on $Z'$, we get an exceptional divisor $E$ with self-intersection $-2$. The proper transform $\til{Z'}$ of $Z'$ has self-intersection $-1$. We can successively blow down $\til{Z'}$ and $E_1$, obtaining a smooth space at each stage. Thus $S_1(4)$ is  smooth.

We have now established the result of the theorem for $S(N)$ and $S_1(N)$ with $N>2$.
For $S(2)$, we need to analyze the fixed points of 
$$
\G(2)/\G(4) = \left\{\pm I, \pm {1\;\;0\choose 2\;\;1}, \pm {1\;\;2\choose 0\;\;1}, 
\pm {1\;\;2\choose 2\;\;1}\right\} 
$$ 
in $T(4)$. First, at any point of $Z$, including cusps, $-I$ 
acts as $(x,t)\to (x,-t)$. 
The other matrices have fixed points only 
at the cusps of $Z$. According to Proposition \ref{fp2}, for each such cusp there is one matrix that acts by $(q,s) \to (-q, s)$ and another (its negative) that acts by  $(q,s) \to (-q, -s)$. The invariants of the action are thus $\{q^2, s^2\}$, so the quotient is smooth. We conclude that $T(2)$ is a smooth space.

As before, let $\pi:T(4)\to T(2)$ be the quotient map, and let $Z' = \pi(Z)$. Then we have:
\begin{eqnarray*}
&&\pi_*Z = 4Z', \;\;\pi^*Z' = 2Z,\\
&&\pi^*(Z'^2) = (\pi^*Z')^2 = 4Z^2 = -8,\\
&&Z'^2 = -1.
\end{eqnarray*}
Thus $Z'$ can be blown down, so $S(2)$ is smooth.

We could carry out the same calculation for $S_1(2)$, but in this case there is a much simpler direct approach. The 2-torsion points of the elliptic curve $y^2=x^3+ax+b$ are just the points where $y=0$. Thus $S_1(2)$ is the subvariety of $\P^2_{x,y}\times \A^2_{a,b}$ given by the equations $(y=0,\; y^2=x^3+ax+b)$. It is obviously smooth. 

\section{Proof of the Main Theorem for $S_0(N)$}
\label{S0N}

As before, to determine whether $S_0(N)$ is smooth at $Q$ we must analyze the singularities of $T_0(N)$ and compute the self-intersection of $Z'$. Note that the only values of $N$ for which $S_0(N)$ could possibly be smooth at $Q$ are those where $Z'\isom X_0(N)$ has genus 0, namely, $N\leq 10$ and for $N=12, 13, 16, 18, 25$.

We begin by computing $Z'^2$. Let $\pi: T(N)\to T_0(N)$ be the quotient map. Then we have:
\begin{eqnarray*}
\deg(\pi)&=& N\phi(N) = N^2\prod_{p|N}\left(1-\frac1{p}\right)\\
\pi_*(Z)&=& \frac{N\phi(N)}{2},\qquad \pi^*(Z') = 2Z \\
Z^2 &=& -\frac{1}{24} N^3\prod_{p|N}\left(1-\frac1{p^2}\right)\\
\pi^*(Z'^2) &=& 4Z^2\\
Z'^2 &=& - \frac{1}6 N\prod_{p|N}\left(1+\frac1{p}\right).
\end{eqnarray*}
Note that $Z'^2<-1$ for $N\geq 6$.

The fixed points of $\G_0(N)$ on $\Hbar$, up to equivalence, are of three types \cite{Sch}:
\begin{itemize}
\item[(i)] {\em Type $\rho = e^{\pi i/3}$}: points fixed by an element of order 6 conjugate to 
$0\;\;-1\choose 1\;\;\;\;1$. By Proposition \ref{fp1}, if $\a\in Z$ corresponds to such a fixed point, then the action on $\T_\a$ is 
$$
(x,t)\to\left(\frac{x}{\r^2}, \frac{t}{\rho}\right).
$$
The invariants are $\{x^3, x^2t^2, xt^4, t^6\}$. The quotient is locally isomorphic to a cone over a twisted cubic. Blowing up at the cone point produces an exceptional divisor with self-intersection $-3$.
\item[(ii)] {\em Type $i$}: points fixed by an element of order 4 conjugate to $0\;\;-1\choose 1\;\;\;\;0$. By Proposition \ref{fp1}, if $\a\in Z$ corresponds to such a fixed point, then the action on $\T_\a$ is 
$$
(x,t)\to(-x, -it).
$$
The invariants are $\{x^2, xt^2, t^4\}$. The quotient has an $A_1$ singularity.

\item[(iii)] {\em Type $\infty$}: cusps on $Z'$ fixed by some element $A\in\G_0(N)$ conjugate to $1\;\;1\choose 0\;\;1$. Note that $-A$ is also in $\G_0(N)$. By Proposition \ref{fp2}, $A$ will act by $(q,s)\to(\z q,s)$ and $-A$ by $(q,s)\to(\z q,-s)$. The invariants of the action are $\{q^N, s^2\}$. The quotient is smooth; thus $T_0(N)$ has no singularities at the cusps of $Z'$. 
\end{itemize}

For a general formula giving the number of fixed points of $\G_0(N)$ of each type, see \cite{Sch}. For our purposes, the following table will suffice. The table summarizes all the numerical information that we have obtained in this section. The first column lists all the values of $N$ that need to be considered. We omit $N=2$, since $S_0(2) = S_1(2)$. We also omit those $N\geq 6$ for which $T_0(N)$ has no singularities, since we already know that $Z'^2<-1$. The second and third columns of the table give the number of fixed points of types $\rho$ and $i$ respectively. The last column gives the self-intersection of the proper transform $\til{Z}'$ of $Z'$ after all the singularities of $T_0(N)$ have been resolved. Fixed points of type $\rho$ and $i$  lead to exceptional divisors with self-intersection $-3$ and $-2$ respectively, and thus contribute $-\frac13$ and $-\frac12$ to $\til{Z}'^2$.
$$
\begin{array}{|c||c|c|c|c|}
\hline
N&	\#\r&		\#i&	Z'^2& 		\til{Z}'^2\\
\hline\hline
3&	1&		0&	-\frac23&		-1\\
\hline
4&	0&		0&	-1&			-1\\
\hline
5&	0&		2&	-1&			-2\\
\hline
7&	2&		0&	-\frac43&		-2\\
\hline
10&	0&		2&	-3&			-4\\
\hline
13&	2&		2&	-\frac73&		-4\\
\hline
25&	0&		2&	-5&			-6\\
\hline
\end{array}
$$
From this table, we can immediately determine the intersection matrix of the curves to be blown down in each $T_0(N)$. $S_0(N)$ will be smooth only if this matrix is unimodular, which happens only for $N=4$. Thus only
$S_0(4)$ and $S_0(2)$ are smooth.

This concludes the proof of the Main Theorem.


\section{Topological perspective}
\label{topological.sec}

From the standpoint of topology, the complex surfaces 
$S(N)$, $S_1(N)$, and $S_0(N)$ naturally give rise to 
covers of $S^3$ branched over the trefoil knot and with
monodromy in $\slzn$.
Our results determine the homeomorphism types of these covers, 
so we now go back and reexamine our problem from this topological
perspective.  We have tried to keep the exposition in this section as
self-contained as possible, but will need to refer to
Sections~\ref{SN} and \ref{S0N} for the main results.

Recall from Section \ref{intro} that $W$ is the family of curves
$y^2=x^3+ax+b$ in ${\C}P^2$, where the curves are parameterized by
points $(a,b)\in\C^2$. Away from the discriminant locus
$\D=\{4a^3+27b^2=0\}$, the fiber over $(a,b)$ is an elliptic
curve. Thus $W$ restricts to a torus bundle $\til{W}$ over
$\C^2\setminus\D$. In Equations~\eqref{SofN.eq}, \eqref{S1ofN.eq}
and~\eqref{S0ofN.eq}, we use $\til{W}$ to define three families of
finite covers of $\C^2\setminus\D$, which we then extend to (possibly
singular) spaces $S(N)$, $S_1(N)$ and $S_0(N)$ mapping surjectively
onto $\C^2$. For simplicity of exposition, we discuss $S_1(N)$ in
detail below and consider $S(N)$ and $S_0(N)$ briefly afterwards.

As an elliptic curve, each fiber $T^2_{(a,b)}$ of $\til{W}$ comes
equipped with a group structure. $\til{S}_1(N)$ consists of the points
of {\em exact} order $N$ in each fiber, which may also be canonically
identified with elements of $H_1(T^2_{(a,b)};\Z/N)$.  The
normalization $S_1(N)$ completes $\til{S}_1(N)$ to a branched cover of
$\C^2\setminus\{(0,0)\}$, with a single point $Q$ above $(0,0)$. Our
goal in the preceding sections was to determine for which $N$ the
varieties $S(N)$, $S_1(N)$ and $S_0(N)$ are smooth over $(0,0)$.

In order to study the topology of $S_1(N)$ at $Q$, we consider the
restriction of this branched cover to the unit sphere $S^3$, with
branch locus $K = S^3\cap\D$, a trefoil knot. Let $\til{M}_1(N)$,
$M_1(N)$ be the restrictions of $S_1(N)$ to
$S^3\setminus K$ and $S^3$ respectively. Since $\C^2\setminus\D$ is
homeomorphic to $(S^3\setminus K)\times \R^+$, a neighborhood of $Q$
in $S_1(N)$ is a cone on $M_1(N)$, and in particular, for $S_1(N)$ to
be smooth at $Q$ it is necessary that $M_1(N)\cong S^3$. However, as
noted in Section~\ref{intro}, the homeomorphism type of $M_1(N)$ (and
not simply whether it is the three-sphere) may be of topological
interest in its own right.

Given a basepoint $(a_0,b_0)\in S^3\setminus K$, the $3$-manifold
$\til{M}_1(N)$ is the cover of $S^3\setminus K$ corresponding to the
stabilizer of $(0,1)\in H_1(T^2_{(a_0,b_0)};\Z/N)$ under the monodromy
action of $\pi_1$.  The fundamental group of the trefoil complement is
the braid group $B_3$, with presentation
\[
B_3 \cong \langle x,y| x^3=y^2 \rangle,
\]
and the homological monodromy of $W$ is known to be
\[
\slz \cong \langle x,y| x^3=y^2, x^6 =1 \rangle
\]
under the obvious map.
For a suitable choice of basis for $H_1$, this is given by
\begin{align*}
x &\mapsto \twobytwo{0}{-1}{1}{1}, & y &\mapsto \twobytwo{0}{-1}{1}{0}
\end{align*}
(note that $\pi_1(S^3\setminus K)$ acts on $H_1(T^2_{(a_0,b_0)})$ on
the right).  Up to conjugacy, the fundamental group of $\til{M}_1(N)$
is therefore the preimage in $B_3$ of the congruence subgroup
\[
\G_1(N) = \left\{ A\in\slz : 
            A\equiv {\textstyle\twobytwo{1}{*}{0}{1}}\bmod N \right\}.
\]

To describe $M_1(N)$, we look at the monodromy of $\til{M}_1(N)$ on
the boundary of a regular neighborhood $U$ of $K$ in $S^3$. As
generators for $\pi_1(\partial U)$ we use a meridian $\mu$ and
longitude $\lambda$ such that $[\mu]=yx^{-1}$ and $[\lambda]=y^2$;
note that $\lambda$ is not the standard choice of longitude, but
rather is a parallel copy of $K$ on the torus on which $K$ lies as a
$(2,3)$-torus knot. Then $\mu$, $\lambda$ act by
$\twobytwo{1}{0}{1}{1}$, $-I$ in $\slzn$ respectively. Thus meridional
discs of $K$ are branch-covered by their preimages in $M_1(N)$, with
branch orders dividing $N$, and $K$ is typically double-covered by the
connected components of its preimage for $N>2$. For example, if $p$ is
an odd prime, then the preimage of the regular neighbourhood $U$ has
$p-1$ connected components in $M_1(p)$. Of these, half map to $S^3$
with no branching, the remaining half map with $p$-fold branching, and
in each component the core circle double-covers $K$.

We now note that level curves of the $j$-invariant
$j=1728\cdot{4a^3}/(4a^3 +27b^2)$ give a Seifert fibration of $S^3$,
isomorphic to that given by the $S^1$ action
$e^{i\theta}\cdot(a,b)=(e^{2i\theta}a,e^{3i\theta}b)$.  The fiber
corresponding to $j=\infty$ is $K$, and there are two exceptional
fibers ($j=0$ and $j=1728$), which parametrize elliptic curves with
automorphisms. The orbit space of $S^3$ is homeomorphic to
$\Hbar/\slz$, where $\Hbar=\H \cup \P^1(\Q)$ is the extended upper
half plane, and $\slz$ acts by fractional linear transformations. The
Seifert fibration lifts to $M_1(N)$, with orbit space homeomorphic to
$X_1(N)=\Hbar/\Gamma_1(N)$.

An exceptional fiber of $M_1(N)$ must be either a preimage of an
exceptional fiber of $S^3$, or a preimage of $K$. The first
possibility occurs when $\G_1(N)$ has torsion, and corresponds to the
fixed point of a torsion element; the second possibility can only
occur when a conjugate of $\lambda\mu^{-k}$ acts by an element of
$\G_1(N)$ for some $0<k<N$.  A matrix $A\in\slz$ has finite order if
and only if $|\tr A| \leq 1$; since elements of $\G_1(N)$ have trace
$2$ mod $N$, there can be torsion elements only for $N=2$ or $3$.
Further, it can be checked that a preimage of $K$ can be an
exceptional fiber only for $N=4$. It follows that for $N\geq 5$,
$M_1(N)$ is a circle bundle, corresponding to the complex line bundle
$\L$ of Section~\ref{SN}.

When $N\geq 5$, the homeomorphism type of $M_1(N)$ is thus completely
determined by the genus of $X_1(N)$ and the 
Chern class of the circle bundle. The genus can be found 
in~\cite[p.\ 219]{DH}, and the Chern class is given in Section~\ref{SN},
Equation~\eqref{c1M1N.eq1}. In contrast, $M_1(2)$, $M_1(3)$, and
$M_1(4)$ have exceptional fibers and must be
analyzed individually. The results of Section~\ref{SN} show that in
all these cases (as well as for $N=5$ and $6$), $M_1(N)$ is
homeomorphic to $S^3$. Each space has a single exceptional 
fiber: those of $M_1(2)$ and $M_1(3)$ are lifts of an exceptional
fiber of $S^3$, and have multiplicities $2$ and $3$ respectively,
while that of $M_1(4)$ is a lift of $K$ and has multiplicity $2$.

In the same way, we can obtain another Seifert fibered space
$M(N)$ branch-covering $S^3$ from the space $\til{S}(N)$ defined in
Equation~\eqref{SofN.eq}.  Points on the fiber of $\til{S}(N)$
over $(a,b)$ correspond to pairs of points in $H_1(T^2_{(a,b)};\Z/N)$
that form a symplectic basis with respect to the intersection pairing, so
the fundamental group of $\til{M}(N)$ is the preimage in $B_3$ 
of the congruence subgroup 
\[
\Gamma(N) = \ker (\slz\rightarrow \slzn).
\]
The orbit space of $M(N)$ is $X(N) = \Hbar/\Gamma(N)$, the genus of
which is known and can be found in \cite[p.\ 219]{DH} or \cite[p.\ 96]{Sch}.

For $N\geq 3$ the Seifert fibration of $M(N)$ is a circle bundle, with
Chern class given by Equation~\eqref{c1MN.eq}. Combined with the genus
of $X(N)$ this determines the homeomorphism type of $M(N)$ for $N\geq 3$, and a
separate calculation in Section~\ref{SN} shows that $M(2)\cong S^3$.

Finally, the space $\til{S}_0(N)$ defined in Equation~\eqref{S0ofN.eq}
yields a third Seifert fibered space, $M_0(N)$. The fiber of
$\til{S}_0(N)$ over $(a,b)$ is the set of cyclic subgroups of
$H_1(T^2_{(a,b)};\Z/N)$ of order $N$, so the fundamental group of
$\til{M}_0(N)$ is the preimage in $B_3$ of the congruence subgroup
\[
\G_0(N) = \left\{ A\in\slz : 
            A\equiv {\textstyle\twobytwo{*}{*}{0}{*}}\bmod N \right\}.
\]
The trace condition no longer holds, so in contrast to $M_1(N)$ and
$M(N)$, $M_0(N)$ has exceptional fibers for infinitely many values of
$N$.  Our results are thus insufficient to determine the homeomorphism
type of $M_0(N)$ in general.  However, Section~\ref{S0N} shows that
$M_0(N) \cong S^3$ if and only if $N=2$ or $4$.


\end{document}